\documentclass[12pt,a4paper]{amsart}
\usepackage{amsmath, amsthm, amssymb, enumerate, mathtools}

\usepackage[breaklinks=true,colorlinks=true,linkcolor=blue,citecolor=red,urlcolor=blue]{hyperref}

\usepackage[margin=1in,footskip=0.25in]{geometry}

\usepackage{color}

\usepackage[utf8]{inputenc}
\usepackage{eepic,epic}
\usepackage{esint}
\usepackage{multirow}
\usepackage{array}
\usepackage[utf8]{inputenc}
\usepackage{graphicx}
\usepackage{caption}
\usepackage{booktabs}
\usepackage{siunitx}
\usepackage{makecell}
\usepackage{mathrsfs}
\usepackage{enumitem}

\usepackage{theoremref}

\def\to{\longrightarrow}
\def\mapsto{\longmapsto}

\def\cwedge{\bigcirc\kern-1.07em\wedge\ }

\newcommand \tM{\overline{M}}

\newcommand \la{\lambda}

\newcommand \<{\langle}
\renewcommand \>{\rangle}

\newcommand \tR{\overline{R}}

\newcommand \Tr{\operatorname{Tr}}

\newcommand \rank{\operatorname{rank}}
\newcommand \id{\operatorname{id}}

\theoremstyle{plain}

\newtheorem{thm}{Theorem}[section]
\newtheorem{lem}[thm]{Lemma}
\newtheorem{prop}[thm]{Proposition}

\usepackage{float}

\theoremstyle{definition}
\newtheorem{defn}{Definition}
\newtheorem*{defn*}{Definition}

\theoremstyle{remark}
\newtheorem{rem}[thm]{Remark}

\numberwithin{equation}{section}

\newcommand{\nb}{\nabla}

\newcommand{\mr}{\mathbb{R}}

\begin{document}
	
	\title[Weakly Einstein submanifolds in space forms]
	{On weakly Einstein submanifolds in space forms satisfying certain equalities}

	\author{Jihun Kim \and JeongHyeong Park}
	\address{Department of Mathematics, Sungkyunkwan University, Suwon, 16419, Korea}
\email{jihunkim@skku.edu, parkj@skku.edu}


\begin{abstract}
We classify weakly Einstein submanifolds in space forms that satisfy Chen's equality. We also give a classification of weakly Einstein hypersurfaces in space forms that satisfy the semisymmetric condition. In addition, we discuss some characterizations of weakly Einstein submanifolds in space forms whose normal connection is flat.
\end{abstract}

	\maketitle

\section{Introduction}
Let $(M,g)$ be an $n$-dimensional Riemannian manifold. We denote the Levi-Civita connection by $\nabla$ and the Riemannian curvature tensor by $R$, where $R(X,Y)=[\nb_X,\nb_Y] - \nb_{[X,Y]}$ for $X,Y\in \mathfrak{X}(M)$ with $\mathfrak{X}(M)$ representing the Lie algebra of all smooth vector fields on $M$. Furthermore, we define the Ricci tensor as $\rho$, given by $\rho(X,Y)=\Tr(Z\mapsto R(Z,X)Y)$, and the scalar curvature as $\tau=\Tr_{g}\rho$. For a local orthonormal frame $\{e_{i}\}$, we define
\begin{gather*}
    \check{R}(X,Y):=\sum\limits_{a,b,c=1}^n R(X,e_{a},e_{b},e_{c})R(Y,e_{a},e_{b},e_{c}),\\
    \check{\rho}(X,Y):=\sum\limits_{a=1}^{n}\rho(X,e_a)\rho(Y,e_a),\quad R[\rho](X,Y):=\sum\limits_{a,b=1}^{n}R(X,e_a,e_b,Y)\rho(e_a,e_b)
\end{gather*}
for vector fields $X$ and $Y$. It is noteworthy that $\check{R}$, $\check{\rho}$, and $R[\rho]$ are a symmetric $(0,2)$-tensor field on $M$. Berger \cite{Ber} showed the following curvature identity holds for a $4$-dimensional compact Riemannian manifold:
\begin{equation}\label{eq:4dim}
    \left(\check{R}-\frac{||R||^2}{4} g\right) -2\left(\check{\rho}-\frac{||\rho||^2}{4}g\right) - 2\left(R[\rho]-\frac{||\rho||^2}{4}g\right) + \tau\left(\rho-\frac{\tau}{4}g\right)=0.
\end{equation}
One of the authors proved that Equation \eqref{eq:4dim} holds on any $4$-dimensional Riemannian manifold \cite{EPS13}, and discussed critical metrics for the squared $L^2$-norm functionals of the curvature tensor \cite{EPS11}. 
Also, the following concept was introduced as a generalization of Einstein manifolds in dimension $4$.
\begin{defn}[\cite{EPS13}]\label{def:wE}
An $n$-dimensional Riemannian manifold $(M,g)$ is said to be \textit{weakly Einstein} if
\begin{equation}\label{eq:wEcondition}
    \check{R}(X,Y) = \frac{||R||^2}{n} g(X,Y)
\end{equation}for all vector fields $X$ and $Y$.
\end{defn}

There is a geometric meaning of the condition~\eqref{eq:wEcondition}: If a compact manifold has a metric with a parallel Ricci tensor (in particular, an Einstein metric), then the metric is critical for the curvature functional $g\mapsto \int_{M}||R||^2 dV_g$ restricted to metrics with $\operatorname{vol}(M)=1$ if and only if the condition \eqref{eq:wEcondition} is satisfied (see \cite[Corollary 4.72]{Bes}).

For $n=3$, it is well-known that $M$ is an Einstein manifold if and only if it has constant curvature. However, $M$ is a weakly Einstein manifold if and only if either it is Einstein or its Ricci operator is of rank one \cite{GHMV}.
In the case $n=4$,  Euh, Park, and Sekigawa (EPS) \cite{EPS14} provided two examples of weakly Einstein manifolds as a direct product of two $2$-dimensional Riemannian manifolds, $M_{1}(c)$ and $M_{2}(-c)$, of constant Gaussian curvature $c$ and $-c$ respectively ($c\neq 0$), along with a space named an ``EPS space"  by Arias-Marcos and Kowalski \cite{AK}. It was later proven that any $4$-dimensional homogeneous weakly Einstein manifold, which is not Einstein, is isometric to one of these two examples \cite{AK}.
We remark that a Riemannian product $M=M_{1}^{n_1}(c_1)\times M_{2}^{n_2}(c_2)$ of two $n_k$-dimensional manifolds of constant curvature $c_k$ ($k=1,2)$ is weakly Einstein
if and only if $c_1^2(n_1-1)=c_2^2(n_2-1)$.
For example, the product manifold $\mathbb{S}^2(\sqrt{2}) \times \mathbb{H}^{3}(-1)$ is weakly Einstein but not Einstein, whereas $\mathbb{S}^2(2)\times \mathbb{S}^3(1)$ is Einstein but not weakly Einstein. Hence, for $n\ge 5$, the Einstein condition and the weakly Einstein condition may not be related to each other. However, it was shown that these two conditions are equivalent for a specific type of manifold known as Golden space form \cite{KPS}.

The weakly Einstein condition has also been used in critical metric problems (see \cite{BDO, HY}). On the other hand, in \cite{GHMV}, the authors classified locally conformally flat weakly Einstein manifolds, and as an application, they obtained the classification of weakly Einstein hypersurfaces in a Euclidean space: a non-Einstein hypersurface in a Euclidean space is weakly Einstein if and only if it is homothetic to an open part of a warped product $I\times_{f}M(c)$, where $I$ is a real interval, $M(c)$ is a space of constant curvature $c$, and $f'(t)^2 + f(t)f''(t)-c=0$. Recently, Wang and Zhang \cite{WZ} classified weakly Einstein real hypersurfaces in a nonflat complex space form of complex dimension two. Specifically, they proved that a $3$-dimensional Hopf hypersurface is weakly Einstein if and only if it is locally congruent to a geodesic sphere with a special radius in the complex hyperbolic plane.

The structure of this paper is as follows: In Section \ref{sec2}, we classify weakly Einstein submanifolds in space forms that satisfy Chen's equality.
In Section \ref{sec4}, we classify weakly Einstein semisymmetric hypersurfaces in space forms.
In Section \ref{sec6}, for a weakly Einstein submanifold in a space form with flat normal connection, we give a sufficient condition to be $2$-stein, and we find all possible forms of the shape operators of a $4$-dimensional weakly Einstein submanifold of codimension two in a Euclidean space.

\section{Weakly Einstein submanifolds satisfying Chen's equality}\label{sec2}
We focus on a submanifold $M^n$ immersed in a space form $\tM^{n+p}(c)$ of constant curvature $c$, where $n\geq 3$ and $p\geq 1$. The inner product $\<\cdot,\cdot\>$ on $\tM^{n+p}(c)$ induces an inner product on $M^n$.
We denote the second fundamental form and the shape operator by $h$ and $A_\xi$ for $\xi\in\Gamma(T^\perp M^n)$, respectively.
Gauss equation is given by
\begin{align*}
R(X,Y,Z,W) &= \tR(X,Y,Z,W) + \<h(X,W),h(Y,Z)\> - \<h(X,Z),h(Y,W)\>\\
&=c\{\<X,W\>\<Y,Z\>-\<X,Z\>\<Y,W\>\} \\
&\quad+ \<h(X,W),h(Y,Z)\> - \<h(X,Z),h(Y,W)\>
\end{align*}
for $X,Y,Z,W\in \Gamma(TM^n)$. Note that $\<h(X,Y),\xi\> = \<A_\xi X,Y\>$.

We define the mean curvature vector as $H = \frac{1}{n}\mathrm{Tr} h$. A submanifold $M^n$ is called \emph{totally geodesic} if $h\equiv 0$, and \emph{minimal} if $H\equiv 0$. For a plane section $\pi = X \wedge Y$ in $T_xM^n$ at a point $x\in M^n$, we define the sectional curvature $K(\pi)$ of $M^n$ as $$K(\pi) = \frac{R(X,Y,Y,X)}{\<X,X\>\<Y,Y\>-\<X,Y\>^2}.$$ We consider the real function $\inf K$ on $M^n$ given by
$$(\inf K)(x) := \inf\{K(\pi):\pi\text{ is a plane in }T_xM^n\}.$$
Chen established an inequality involving the intrinsic scalar invariants $\inf K$ and $\tau$ and the extrinsic scalar invariant $||H||$ as follows:
\begin{lem}[\cite{Chen}]\label{lem:Chen}  
    Let $M^n$, $n\ge 2$, be a submanifold of a space form $\tM^{n+p}(c)$, $p\ge 1$. Then
\begin{equation}\label{eq:Chen}
    \frac{1}{2}\tau -\inf K   \le  \frac{n^2(n-2)}{2(n-1)}||H||^2 + (n+1)(n-2)\frac{c}{2}.
\end{equation}
Equality holds at a point $x$ if and only if with respect to suitable local orthonormal frames $\{e_i\}_{i=1,\ldots,n}$ for $T_x M^n$ and $\{\xi_j = e_{n+j}\}_{j=1,\ldots,p}$ for $T_x^\perp M^n$, the shape operator $A_t:=A_{\xi_t}$ are given by
\begin{equation}\label{eq:Cheneq}
    A_1 =
    \begin{pmatrix}
    a & 0 & 0 & 0 &\cdots &0\\
    0 & b & 0 & 0 &\cdots &0\\
    0 & 0 & \mu & 0 &\cdots &0\\
    0 & 0 & 0 & \mu &\cdots &0\\
    \vdots & \vdots & \vdots & \vdots &\ddots &0\\
    0 & 0 & 0 & 0 &\cdots &\mu\\
    \end{pmatrix},
    \quad A_t  =
    \begin{pmatrix}
    c_t & d_t & 0 & 0 &\cdots &0\\
    d_t & -c_t & 0 & 0 &\cdots &0\\
    0 & 0 & 0 & 0 &\cdots &0\\
    0 & 0 & 0 & 0 &\cdots &0\\
    \vdots & \vdots & \vdots & \vdots &\ddots &0\\
    0 & 0 & 0 & 0 &\cdots &0\\
    \end{pmatrix},
    (t\ge 2),
\end{equation}
where $\mu = a+b$. For any such frame, $\inf K(x)$ is attained by the plane $e_1 \wedge e_2$.

The inequality \eqref{eq:Chen} is known as Chen's (basic) inequality. When equality holds, it is called Chen's (basic) equality {\rm (}and if the equality holds at any point $x \in M^n$, the submanifold $M^n$ is called Chen $\delta(2)$-ideal submanifold \cite{Chen2}{\rm )}.
For $n=2$, Chen's equality trivially holds.
\end{lem}

\begin{rem}\label{rem:formula}
Let $c=0$, that is, the ambient space is a Euclidean space. From Gauss equation, we have the following formulas for $K_{rs} := K(e_r \wedge e_s)$ (assuming the submanifold to be Chen $\delta(2)$-ideal):
\begin{gather*}
    K_{12} = ab-\sum\limits_{t=2}^{p}(c_t^2 + d_t^2),\quad K_{1j} = a\mu,\quad K_{2j} = b\mu,\quad K_{ij} = \mu^2\;(i\neq j),\\
    \rho_{11} = K_{12} + (n-2)a\mu,\quad\rho_{22}=K_{12}+(n-2)b\mu,\quad\rho_{ii}=(n-2)\mu^2,\quad \rho_{ij}=0\;(i\neq j)\\
    \tau = 2K_{12} + (n-1)(n-2)\mu^2.
\end{gather*}
where $3 \le i,j \le n$. Furthermore, $R(e_i,e_j)e_k = 0$ if $i,j$ and $k$ are mutually different \cite{DPV}.
\end{rem}

As a result, Dillen et al. proved the following:
\begin{thm}[\cite{DPV}]
    Let $M^n$ be a submanifold in $\mr^{n+p}$ satisfying Chen's equality. Then $M^n$ is Einstein if and only if $M^n$ is totally geodesic.
\end{thm}

In this paper, we consider the weakly Einstein condition~\eqref{eq:wEcondition},
for a submanifold of a space form that satisfies Chen's equality. We will use the same notation for $a,b,\mu,c_t$, and $d_t$ as in \eqref{eq:Cheneq}.
\begin{thm}\label{thm:cheneuclid}
    Let $M^n$ be a weakly Einstein submanifold in a Euclidean space $\mr^{n+p}$ satisfying Chen's equality. Then, if $p=1$, the hypersurface $M^n$ is flat.
    If $p\ge 2$, then either
    \begin{enumerate}[label={\rm{(\roman*)}}]
        \item $M^n$ is flat, or
        \item $a=b\neq 0$ and the Ricci curvature is non-negative. If $n=4$, we additionally obtain that $a= b=\pm\frac{1}{2}$ and the Ricci operator has rank $2$.
    \end{enumerate}
\begin{proof}
     Since $R(e_i,e_j)e_k = 0$ if $i,j$ and $k$ are mutually different, we get
    \begin{align*}
        &\check{R}_{11} = K_{12}^2 + K_{13}^2 + \cdots + K_{1n}^2 = K_{12}^2 + a^2\mu^2(n-2) \\
        &\check{R}_{22} = K_{21}^2 + K_{23}^2 + \cdots + K_{2n}^2 = K_{12}^2 + b^2\mu^2(n-2)  \\
        &\check{R}_{33} = K_{31}^2 + K_{32}^2 + \cdots + K_{3n}^2  = a^2\mu^2 + b^2\mu^2 + \mu^4(n-3) \\
        &\qquad\vdots\\
        &\check{R}_{nn} = K_{n1}^2 + K_{n2}^2 + \cdots + K_{n n-1}^2 = a^2\mu^2 + b^2\mu^2 + \mu^4(n-3)\; (n\ge 4).
    \end{align*}
    Assuming $M^n$ is weakly Einstein, we obtain
    \begin{equation}\label{eq:wEchen}
        K_{12}^2 + a^2\mu^2(n-2) =  K_{12}^2 + b^2\mu^2(n-2) = a^2\mu^2 + b^2\mu^2 + \mu^4(n-3).
    \end{equation}
    From the first equality of \eqref{eq:wEchen}, we get $(a+b)(a-b)\mu=\mu^2(a-b)=0$, which yields either $\mu=0$ or $a=b$. When $\mu=0$, we obtain $K_{ij}=0$ for $1\le i,j\le n$, and therefore $K\equiv 0$ and $M^n$ is flat.

      Now we suppose that $\mu\neq 0$. For $a=b\neq0$, using the second equality of \eqref{eq:wEchen}, we get $K_{12}^2 = 4a^4 (3n-8)$ and hence $K_{12}=\pm 2a^2\sqrt{3n-8}$. We note that if $p=1$, any of them can not hold because $K_{12}=a^2$. If $K_{12}=2a^2\sqrt{3n-8}$, then $\sum\limits_{t=2}^{p}(c_t^2 + d_t^2)= a^2 (1-2\sqrt{3n-8})<0$, a contradiction. Hence we have $K_{12} = -2a^2\sqrt{3n-8}$. The Ricci curvatures are then given by
    \begin{equation*}
        \rho_{11}=\rho_{22}=2a^2(n-2 - \sqrt{3n-8}),\;\rho_{ii}=4(n-2)a^2>0\;(i\ge3).
    \end{equation*}
    Here, $n-2-\sqrt{3n-8}\ge0$ (resp. $<0$) if and only if $n\le3$ or $n\ge4$ (resp. $3<n<4$). Therefore, there does not exist negative Ricci curvature. Furthermore, we observe that  if $n$ equals either $3$ or $4$, $\rho_{11}=\rho_{22}=0$, with the rank of the Ricci operator being $1$ and $2$, respectively. Assuming $n=4$, we have $K_{12}^2 = 16a^4$.  In \cite{EPS14}, the authors established the existence of a generalized Singer-Thorpe basis such that
    \begin{equation}\label{eq:gst}
        R_{1212}^2 = R_{3434}^2,\;R_{1313}^2 = R_{2424}^2,\;R_{1414}^2 = R_{2323}^2,
    \end{equation}
    and $R_{ijjk}=0$ for $i\neq k$, $1\le i,j,k\le 4$ (see Lemma~\ref{lem:gST}). Using the relation $K_{12}^2 =R_{1212}^2 = R_{3434}^2 = K_{34}^2$ and $K_{34}^2 = 4a^2$, we obtain $16a^4 = 4a^2$, from which it follows that $a^2= \frac{1}{4}$.
    \end{proof}
\end{thm}

\medskip

We now study a submanifold of a nonflat space form $\tM^{n+p}(c)$ satisfying Chen's equality.
Similar to Remark \ref{rem:formula}, we have the following formulas for $K_{rs}$: 
\begin{gather*}
    K_{12} = c+ab-\sum\limits_{t=2}^{p}(c_t^2 + d_t^2),\; K_{1j} = c+ a\mu,\; K_{2j} =c+ b\mu,\; K_{ij} = c+\mu^2,\\
    \rho_{11} = K_{12} + (n-2)(c+a\mu),\quad\rho_{22}=K_{12}+(n-2)(c+b\mu),\quad\rho_{ii}=c+(n-2)(c+\mu^2),\\
    \tau = 2K_{12} + (n-2)(2c+(n-1)(c+\mu^2)),
\end{gather*}
where $3 \le i,j \le n$. Furthermore, $R(e_i,e_j)e_k = 0$ if $i,j$ and $k$ are mutually different.
Also, we note that if $\mu = 0$, then $\Tr(A_1)=\Tr(A_t)=0$ $(t\ge 2)$ and so $M^n$ is minimal.

\begin{thm}\label{thm:chenspaceform}
Let $M^n$ be a weakly Einstein submanifold in a nonflat space form $\tM^{n+p}(c)$ satisfying Chen's equality.

When $p=1$ (i.e. when $M^n$ is a hypersurface), we have
    \begin{itemize}
        \item[{\rm{(i)}}] if $c>0$, then either
        \begin{itemize}
            \item[{\rm{(i.a)}}] $M^n$ is a space of constant curvature $c$, or
            \item[{\rm{(i.b)}}] $M^n$ is minimal of non-negative Ricci curvature, and $(a,b)=(\pm\sqrt{2c},\mp\sqrt{2c})$;
        \end{itemize}
        \item[{\rm{(ii)}}] if $c<0$, then either
        \begin{itemize}
           \item[{\rm{(ii.a)}}] $M^n$ is a space of constant curvature $c$, or
           \item[{\rm{(ii.b)}}] either $n=3$ or $n=4$ and then $a=b=\pm\sqrt{-\frac{c(4n-10)}{12n-33}}$, or
            \item[{\rm{(ii.c)}}] $b=-a\pm\sqrt{-2c}$.
        \end{itemize}
    \end{itemize}
When $p\ge 2$, we have
    \begin{itemize}
        \item[{\rm{(iii)}}] if $c>0$, then either
        \begin{itemize}
            \item[{\rm{(iii.a)}}] $M^n$ is a space of constant curvature $c$, or
            \item[{\rm{(iii.b)}}] $M^n$ is minimal of non-negative Ricci curvature, or
            \item[{\rm{(iii.c)}}] either $n=3$ or $n=4$;
        \end{itemize}
        \item[{\rm{(iv)}}] if $c<0$, then either
        \begin{itemize}
        \item[{\rm{(iv.a)}}] $M^n$ is a space of constant curvature $c$, or
            \item[{\rm{(iv.b)}}] either $n=3$ or $n=4$, or
            \item[{\rm{(iv.c)}}] $b=-a\pm\sqrt{-2c}$.
        \end{itemize}
    \end{itemize}
\begin{proof}
By applying the weakly Einstein condition, we derive the following equations:
     \begin{equation}\label{eq:spwEchen}
        K_{12}^2 + (c+a\mu)^2(n-2) =  K_{12}^2 + (c+b\mu)^2(n-2) = (c+a\mu)^2 + (c+b\mu)^2 + (c+\mu^2)^2(n-3).
    \end{equation}
    As a consequence of the first equality of~\eqref{eq:spwEchen}, we can deduce that $c+a\mu=\pm(c+b\mu)$, which implies $(a-b)\mu=0$ or $\mu^2=-2c$.

    \smallskip
    \noindent\textsc{Case 1.} Let $c>0$. Then we have either $\mu=0$ or $a=b$. If $\mu=0$, then $K_{12}^2=c^2$ and so we get $K_{12}=\pm c$. If $K_{12}=c$, then $K_{ij}=c$ for $1\le i,j\le n$ and $M^n$ is a space of constant curvature $c$. Suppose $K_{12}=-c$, then we get
    \begin{equation*}
        \rho_{11}=\rho_{22}=c(n-3),\quad \rho_{ii}=c(n-1)\;(i\ge 3).
    \end{equation*}
    Thus, the Ricci curvature is non-negative. When $p=1$, since $K_{12}=c-a^2$, we find that $a^2=2c$ if $K_{12}=-c$.

   We now assume that $\mu\neq 0$. For $a=b\neq0$, using the second equality of~\eqref{eq:spwEchen} we obtain
    $$K_{12}^2 = (12n-32)a^4 + (4n-8)ca^2+c^2,$$
    and since $K_{12}^2\ge 0$ for $a^2$, the discriminant $D=(2(n-2)c)^2-c^2(12n-32)\le 0$, which implies $3\le n \le 4$. Thus $n=3$ or $n=4$. When $p=1$, by substituting $K_{12}=c+a^2$, we have $a^2 = -\frac{c(4n-10)}{12n-33}$. But this does not occur for $c>0$.

    \smallskip
    \noindent\textsc{Case 2.} Let us consider the case where $c<0$. If $\mu=0$, the situation is analogous to the case where $c>0$. However, if $K_{12}=-c$, we get a contradiction because $0<-2c=-a^2 - \sum\limits_{t=2}^{p}(c_t^2 + d_t^2)$. If $a=b$, we obtain a result similar to that of the $c>0$ case. Finally, if $\mu^2 = -2c$, we have $b=-a\pm\sqrt{-2c}$.
    \end{proof}
\end{thm}

\section{Weakly Einstein semisymmetric hypersurfaces} \label{sec4}

In this section, we focus on a hypersurface $M^n$ of a space form $\tM^{n+1}(c)$, where $n\ge 3$. Consider the shape operator $A$ with respect to an orthonormal basis in the normal space. By using Gauss equation, we obtain
\begin{equation*}
    R(X,e_{a},e_{b},e_{c}) = c(\<e_{a},e_{b}\>\<X,e_{c}\> - \<e_{a},e_{c}\>\<X,e_{b}\>)
    +\big(\<A e_{a},e_{b}\>\<A X,e_{c}\>-\<A e_{a},e_{c}\>\<A X, e_{b}\>\big)
\end{equation*}
and
\begin{equation}\label{eq:Rcheck-p}
        \check{R}(X,Y) = 2c^2 (n-1)\<X,Y\>+4c\<\big(\Tr(A) A - A^2\big) X,Y\>+2\<\big(\Tr(A^2)A^2 - A^4 \big) X,Y\>.
\end{equation}
Let $\kappa_i$ and $e_i$, for $i=1,\ldots,n$, be the principal curvatures and corresponding principal vectors of $A$, respectively. We substitute $X=Y=e_i$ into~\eqref{eq:Rcheck-p}. As a result, we get
\begin{equation*}
    \check{R}(e_i,e_i)=2c^2 (n-1) + 4c(\Tr(A)\kappa_i-\kappa_i^2)+2(\Tr(A^2)\kappa_i^2-\kappa_i^4).
\end{equation*}

From the weakly Einstein conditon $\check{R}(e_i,e_i)=\check{R}(e_j,e_j)$, we have
\begin{equation}\label{eq:wEceigen}
    (\kappa_i-\kappa_j)\{2c\Tr(A)+(\kappa_i+\kappa_j)(-2c+\Tr(A^2)-\kappa_i^2 -\kappa_j^2\}=0.
\end{equation}
In particular, when $c=0$,
\begin{equation}\label{eq:wEeigen}
    0 = (\kappa_i - \kappa_j )(\kappa_i + \kappa_j )(\Tr(A^2) - \kappa_i^2 - \kappa_j^2),
\end{equation}
for $i,j=1,\ldots,n$.

\smallskip

The definition of a weakly Einstein manifold in \cite{GHMV} is different from our Definition \ref{def:wE} (the authors of \cite{GHMV} additionally require a weakly Einstein manifold to be non-Einstein), 
but we adopt the original definition given in \cite{EPS13} which is our Definition \ref{def:wE}, so
we have a modified version of Lemmas from \cite{GHMV}.

\begin{lem}\label{lem:hypwEEu}
Let $M^n$ be a hypersurface of $\mr^{n+1}$. Then $M^n$ is weakly Einstein if and only if either it is totally umbilical, or flat, or $A$ has exactly two distinct principal curvatures $\pm\kappa$.
\begin{proof}
We can derive from Equation \eqref{eq:wEeigen} that either $M^n$ is totally umbilical, or $A$ has either rank one or exactly two distinct principal curvatures with opposite signs. Note that if $A$ has rank one, then $M^n$ must be flat.
\end{proof}
\end{lem}

\begin{lem}\label{lem:hypwEc}
Suppose $M^n$ is a non-Einstein weakly Einstein hypersurface of $\tM^{n+1}(c)$ with $c\neq 0$. If $A$ has exactly two distinct principal curvatures, then both principal curvatures have the multiplicity greater than one.
\end{lem}

\begin{defn}[\cite{Sz1}]
A Riemannian manifold $(M,g)$ is said to be \textit{semisymmetric} if
\begin{equation}\label{eq:semisym}
    R(X,Y)\cdot R = 0
\end{equation}
for all tangent vector fields $X$ and $Y$, where the endomorphism $R(X,Y)$ acts as a derivation on $R$, that is, for a pair of vector fields $U$ and $V$ we have
\begin{equation*}
    \big(R(X,Y)\cdot R)(U,V) = [R(X,Y),R(U,V)]-R\big(R(X,Y)U,V\big) - R\big(U,R(X,Y)V\big).
\end{equation*}
\end{defn}
The term “semi-symmetric” comes from the fact that the curvature tensor $R_x$ of $M$ at a point
$x \in M$ is the same as the curvature tensor of a symmetric space \cite{Cal}.

We now examine the consequences of the semisymmetric condition
for hypersurfaces $M^n$ in space forms $\tM^{n+1}(c)$. From Gauss equation,
\begin{equation*}
    R(X,Y) = c X\wedge Y + AX \wedge AY,
\end{equation*}
where $X\wedge Y$ denotes the endomorphism which maps $Z$ upon $\<Y,Z\>X - \<X,Z\>Y$.
Then we obtain
\begin{equation*}
    R(e_{i},e_{j}) = (\kappa_{i}\kappa_{j}+c)e_{i}\wedge e_{j}.
\end{equation*}
It follows that
\begin{align*}
    \big(R(e_{i},e_{j})\cdot R)(e_{k},e_{l}) &= [R(e_{i},e_{j}),R(e_{k},e_{l})]\\
    &\quad -R\big(R(e_{i},e_{j})e_{k},e_{l}\big) - R\big(e_{k},R(e_{i},e_{j})e_{l}\big),
\end{align*}
which is equal to zero, except for $k=i$ and $l \neq i,j$. In that case, we get
\begin{align*}
    \big(R(e_{i},e_{j})\cdot R\big)(e_{i},e_{l}) &= (\kappa_{i}\kappa_{j}+c)(\kappa_{i}\kappa_{l}+c)e_{l}\wedge e_{j}\\
    &\quad -(\kappa_{i}\kappa_{j}+c)(\kappa_{j}\kappa_{l}+c)e_{l}\wedge e_{j}.
\end{align*}
Thus, the semisymmetric condition is equivalent to
\begin{equation*}
    (\kappa_{i}\kappa_{j}+c)(\kappa_{i}-\kappa_{j})\kappa_{l}=0,\; {\text{where}}\; 1 \le i, j, l \le n\;  {\text{are mutually different.}}
\end{equation*}

Consequently, the following results have been obtained in \cite{Ry1}.

\begin{lem}\label{lem:hypsemieig}
Let $M^n$ be a hypersurface of a space form $\tM^{n+1}(c)$. Then $M^n$ is semisymmetric if and only if the principal curvatures $\{\kappa_{q}\}_{q=1,\ldots,n}  $ of the shape operator $A$ satisfies $(\kappa_{i}\kappa_{j}+c)(\kappa_{i}-\kappa_{j})\kappa_{l}=0$ for distinct $i,j,l$.
\end{lem}

\begin{prop}\label{prop:semihyprank}
Let $M^n$ be a hypersurface of a space form $\tM^{n+1}(c)$ with $c\neq 0$, and $M^n$ satisfies the semisymmetric condition \eqref{eq:semisym}. Then for any point in $M^n$, either $\rank(A) = n$ or $\rank(A)\le 1$. Moreover, if $\rank(A)=n$, at most two principal curvatures are distinct.

\end{prop}

Ryan \cite{Ry1} also showed the following.
\begin{prop}
    Let $M^n$, $n\ge3$, be a hypersurface in a space form $\tM^{n+1}(c)$ with $c\neq 0$, and $M^n$ satisfies the semisymmetric condition \eqref{eq:semisym}. If $M^n$ has exactly two principal curvatures with multiplicities greater than one, then $M^n$ is locally isometric to a product of two spaces of constant curvature.
\end{prop}

By considering the weakly Einstein condition, we obtain a stronger classification for hypersurfaces as follows. Since the classification of Einstein hypersurfaces in spaces of constant curvature is well-known \cite{Fia}, we now concentrate on non-Einstein hypersurfaces.

\begin{thm}\label{thmA}
Let $M^n$ be a non-Einstein weakly Einstein hypersurface of a space form $\tM^{n+1}(c)$. Suppose that $M^n$ satisfies the semisymmetric condition \eqref{eq:semisym}. Then $c \neq 0$, and we have the following:
    \begin{enumerate}[label={\rm(\roman*)}]
        \item If $c>0$, then $M^n$ is locally isometric to  $\mathbb{S}^{p}(\sin^{-2}\theta) \times \mathbb{S}^{q}(\cos^{-2}\theta)$, where $\cos\theta$, $\sin\theta>0$ and $\tan^4\theta = \frac{p-1}{q-1}$.
        \item If $c<0$, then $M^n$ is locally isometric to  $\mathbb{S}^{p}(\sinh^{-2}\theta)\times \mathbb{H}^{q}(-\cosh^{-2}\theta)$, where $\sinh\theta>0$ and $\tanh^4\theta = \frac{p-1}{q-1}$.
    \end{enumerate}
    In both cases, $n=p+q$ and $p\neq q$ (so $n$ can not be even).
\begin{proof}
The proof presented below is similar to that of \cite{GHMV}.
Assume that $c=0$.
According to Lemma~\ref{lem:hypwEEu}, $A$ has principal curvatures $\pm \kappa$. Thus, each $\kappa_{i}$ in Lemma~\ref{lem:hypsemieig} is either $\kappa$ or $-\kappa$, which implies that $\kappa=0$. Therefore, $A \equiv 0$ and so $M^n$ is flat, a contradiction.

Now let us assume that $c \neq 0$. From Lemma~\ref{lem:hypwEc} and Proposition~\ref{prop:semihyprank}, we have either $A\equiv 0$, or $A = \kappa \id (\neq 0)$, or $A$ has exactly two distinct non-zero principal curvatures of multiplicity greater than one. If $A \equiv 0$ or $A=\kappa\id(\neq 0)$, then $M^n$ has constant curvature and is therefore Einstein.
For the last case, the weakly Einstein condition implies that the two distinct principal curvatures are related. Additionally, as $A$ is a Codazzi tensor, the principal curvatures are constant \cite{Bes,Mer}. Therefore, $M^n$ is an isoparametric hypersurface.
According to the results of isoparametric hypersurfaces in spheres and hyperbolic spaces \cite{CR}, for $n=p+q$, $M^n$ is locally isometric to either $\mathbb{S}^{p}(\sin^{-2}\theta) \times \mathbb{S}^{q}(\cos^{-2}\theta)$ with $\cos\theta,\sin\theta>0$ when $c>0$, or to $\mathbb{S}^{p}(\sinh^{-2}\theta)\times \mathbb{H}^{q}(-\cosh^{-2}\theta)$ with $\sinh\theta>0$ when $c<0$.
As $\mathbb{S}^{p}(\sin^{-2}\theta) \times \mathbb{S}^{q}(\cos^{-2}\theta)$ (resp. $\mathbb{S}^{p}(\sinh^{-2}\theta)\times \mathbb{H}^{q}(-\cosh^{-2}\theta)$) is weakly Einstein,  we have $\tan^{4}(\theta) = \frac{p-1}{q-1}$ (resp. $\tanh^{4}(\theta) = \frac{p-1}{q-1}$).
If $p=q$, in the former case, we have $\cos\theta = \sin\theta$, which means that $M^n$ is Einstein. In the latter case, as $\theta \rightarrow \infty$, the sectional curvatures vanish, which is impossible. Thus, this completes the proof.
\end{proof}
\end{thm}


\section{Weakly Einstein submanifolds with flat normal connection}\label{sec6}

We now consider a Riemannian submanifold $M^n$ in a space form $\tM^{n+p}(c)$ where $n\ge 3$.
For $x \in M^n$ and $X \in T_xM^n$, we use $R_X$ to denote the Jacobi operator defined by $R_XY=R(Y,X)X$ for $Y \in T_xM^n$. A Riemannian manifold $M^n$ is called \emph{$2$-stein} if there exist two functions $f_1, f_2$ on $M^n$ such that
 \begin{align}
     &\Tr R_X = f_1 \|X\|^2,\label{eq:Einstein}\\
    &\Tr (R_X^2) = f_2\|X\|^4,\label{eq:2-stein}
 \end{align}
 for all $x \in M^n$ and all $X \in T_xM^n$ (see \cite{Bes2, CGW}).
 We remark that $f_1$ is a constant when $n \ge 3$, and $f_2$ is a constant when $n \ge 5$ (see  \cite[6.61]{Bes2})

Let $A_t$ be the shape operators relative to an orthonormal basis $\{\xi_t\}_{t=1, \dots, p} $ in the normal space at a point $x \in M^n$. From Gauss equation we have
\begin{equation}\label{eq:RX}
  R_X = c (\|X\|^2 \id - X \otimes X^\flat) + S_X,
\end{equation}
where $S_X$ is given by
\begin{equation}\label{eq:RX'}
    S_X = \sum_{t=1}^{p} (\<A_t X, X\> A_t - (A_t X) \otimes (A_t X)^\flat).
\end{equation}
Therefore, we see that $M^n$ is $2$-stein if and only if at every point $x \in M^n$,
\begin{equation}\label{eq:2stdash}
  \Tr S_X = h_1 \|X\|^2, \qquad \Tr (S_X^2) = h_2 \|X\|^4,  \qquad
\end{equation}
 for some $h_1, h_2 \in \mathbb{R}$ and for all $X \in T_xM^n$.
We substitute the expression for $S_X$ from \eqref{eq:RX'} into \eqref{eq:2stdash} and then we obtain the following.
\begin{prop}[\cite{EKNP}]
    Let $M^n$ be a submanifold in a space form $\tM^{n+p}(c)$. Then $M^n$ is $2$-stein if and only if for any point $x\in M^n$, there exist $h_1, h_2 \in \mathbb{R}$ such that the shape operators $A_t$, $t=1,\ldots,p$, satisfy
    \begin{gather}\label{eq:1stA}
  \sum_{t=1}^{p} \<(\Tr(A_t) A_t - (A_t)^2)X,X\> = h_1 \|X\|^2, \\
  \begin{multlined}[t]
  \sum_{t,s=1}^{p} (\<A_t X, X\> \<A_s X, X\> \Tr(A_t A_s) + \<A_t X, A_s X\>^2 \\
  - \<A_t X, X\> \<A_s A_t A_s X, X\> - \<A_s X, X\> \<A_t A_s A_t X, X\>) = h_2 \|X\|^4,  \label{eq:2stA}
  \end{multlined}
\end{gather}
for all $X \in T_xM^n$.
\end{prop}
We also have
\begin{equation}\label{eq:Rcheck}
    \begin{aligned}
        \check{R}(X,X) &= 2c^2 (n-1)\<X,X\>\\
        &\quad +4c\sum_{t=1}^{p}\<\big(\Tr(A_t) A_t - (A_t)^2\big) X,X\>\\
        &\quad +2\sum_{t,s=1}^{p}\<\big(\Tr(A_t A_s)A_t A_s - A_t A_s A_t A_s\big) X,X\>.
    \end{aligned}
\end{equation}

\begin{thm}\label{thm:suff2stein}
    Let $M^n$ be a weakly Einstein submanifold in a space form $\tM^{n+p}(c)$ whose normal connection is flat.
    \begin{enumerate}[label={\rm{(\roman*)}}]
        \item If $M^n$ is Einstein, then 
         $M^n$ 
         is $2$-stein.
        \item If $M^n$ satisfies the condition~\eqref{eq:2-stein} and $c\neq 0$, then 
        $M^n$ is $2$-stein.
    \end{enumerate}
\begin{proof}
Since the submanifold $M^n$ has flat normal connection, all the shape operators at any point $x \in M^n$ are simultaneously diagonalizable relative to some orthonormal basis $e_i, \; i=1, \dots, n$, for $T_xM^n$. We denote $\la_{t,i}:=\<A_t e_i, e_i\>$. Then equations~\eqref{eq:1stA}, \eqref{eq:2stA}, and \eqref{eq:Rcheck} are equivalent to
\begin{gather}\label{eq:flEin}
  \sum_{t=1}^{p} (\Tr(A_t) \la_{t,i} - (\la_{t,i})^2) = h_1, \\
  \sum_{t,s=1}^{p} (\Tr(A_t A_s) \la_{t,i} \la_{s,i} - (\la_{t,i} \la_{s,i})^2) = h_2, \label{eq:fl2st}\\
  \check{R}_{ii} = 2c^2 (n-1) + 4c \sum_{t=1}^{p} (\Tr(A_t) \la_{t,i} - (\la_{t,i})^2) + 2 \sum_{t,s=1}^{p} (\Tr(A_t A_s) \la_{t,i} \la_{s,i} - (\la_{t,i} \la_{s,i})^2),\label{eq:flweak}
\end{gather}
respectively, for all $i=1,\ldots,n$. We assume that $M^n$ satisfies \eqref{eq:Einstein}. Then from~\eqref{eq:flEin} and \eqref{eq:flweak}, we obtain
    \begin{equation*}
        \sum_{t,s=1}^{p} (\Tr(A_t A_s) \la_{t,i} \la_{s,i} - (\la_{t,i} \la_{s,i})^2) = h_2,
    \end{equation*}
    where $h_2 = \frac{||R||^2}{2n}-c^2(n-1)-2c\, h_1$. Hence, \eqref{eq:2stA} holds and so $M^n$ is $2$-stein.

    Now we assume that $M^n$ satisfies~\eqref{eq:2-stein} and $c\neq 0$. Then from~\eqref{eq:fl2st} and \eqref{eq:flweak}, we get
    \begin{equation*}
        \sum_{t=1}^{p} (\Tr(A_t) \la_{t,i} - (\la_{t,i})^2) = h_1,
    \end{equation*}
    where $h_1 = \frac{1}{4c}\left\{\frac{||R||^2}{n} - 2c^2(n-1) - 2h_2\right\}$. This implies that \eqref{eq:1stA} holds.
    Thus $M^n$ is $2$-stein.
\end{proof}
\end{thm}

\begin{rem}
        In \cite{EKNP}, it was shown that a $2$-stein submanifold in a space form with flat normal connection has constant sectional curvature. Therefore, Theorem~\ref{thm:suff2stein} implies that $M^n$ has constant sectional curvature. 
\end{rem}

\smallskip

We now focus on a special case of $n=4$, $p=2$, and $c=0$.
In \cite{VZ}, the authors obtained all possible forms of the shape operators for a $4$-dimensional Einstein submanifold in $\mathbb{R}^6$ whose normal connection is flat.
In this paper, we study a $4$-dimensional weakly Einstein submanifold in a Euclidean space $\mathbb{R}^6$ whose normal connection is flat.

We recall that a submanifold $M^n$ is called \emph{quasi-umbilical} in the direction $\xi$ if the shape operator $A_\xi$ admits an eigenvalue $\la$ with the multiplicity $n-1$ or $n$. If $\la=0$, $M^n$ is called \emph{cylindrical} in the direction $\xi$, and $M^n$ is called (\emph{totally}) \emph{cylindrical} if at each point, there exists an orthonormal normal frame with cylindrical directions.

 \begin{lem}[\cite{EPS14}]\label{lem:gST}
Let $M$ be a $4$-dimensional Riemannian manifold. Then $M$ is weakly Einstein if and only if there exists a generalized Singer-Thorpe basis of $T_x M$ at each point $x\in M$ such that
\begin{equation}\label{eq:weak}
    \begin{gathered}
    R_{1212}=a,\quad R_{1313}=b,\quad  R_{1414}=c\\
    R_{3434}=a',\quad  R_{2424}=b',\quad  R_{2323}=c'\\
    R_{1234}=d,\quad  R_{1342}=e,\quad  R_{1423}=-(d+e),
\end{gathered}
\end{equation}
otherwise zero and the relation $a^2 = (a')^2, b^2 = (b')^2, c^2 = (c')^2$ holds. Moreover, if $a=a',b=b',c=c'$, then $M$ is Einstein.
\end{lem}

We prove the following.
\begin{thm}
    Let $M$ be a $4$-dimensional submanifold in $\mathbb{R}^6$ with flat normal connection. Then $M$ is weakly Einstein if and only if for each point $x\in M$, either
\begin{enumerate}[label={\rm{(\roman*)}}]
    \item $M$ is cylindrical, or
    \item $M$ is umbilical with respect to a normal vector $\xi_1$ and cylindrical in $\xi_2$ perpendicular to $\xi_1$, or
    \item for a suitable orthonormal tangent frame of $M$ and an orthonormal normal frame $\{\xi_1,\xi_2\}$, the shape operators $A_1$, $A_2$ are one of the following forms:
    \begin{equation}\label{eq:shape1}
        A_1 =
        \begin{pmatrix}
            \alpha & 0 & 0 & 0\\
            0 & \beta & 0 & 0\\
            0 & 0 & 0 & 0\\
            0 & 0 & 0 & 0
        \end{pmatrix},\quad
        A_2 =
        \begin{pmatrix}
            0 & 0 & 0 & 0\\
            0 & 0  & 0 & 0\\
            0 & 0 & p & 0\\
            0 & 0 & 0 & q
        \end{pmatrix},
    \end{equation}
    where $pq = \pm\alpha\beta\neq 0$;
    \begin{equation}\label{eq:shape2}
        A_1 =
        \begin{pmatrix}
            \alpha & 0 & 0 & 0\\
            0 & 0 & 0 & 0\\
            0 & 0 & \frac{\beta}{\alpha} & 0\\
            0 & 0 & 0 & \frac{\gamma}{\alpha}
        \end{pmatrix},\quad
        A_2 =
        \begin{pmatrix}
            0 & 0 & 0 & 0\\
            0 & \pm\alpha  & 0 & 0\\
            0 & 0 & p & 0\\
            0 & 0 & 0 & q
        \end{pmatrix},
    \end{equation}
    where $pq = -\frac{\beta\gamma}{\alpha^2}\neq 0$;
     \begin{equation}\label{eq:shape3}
        A_1=
        \begin{pmatrix}
            \alpha & 0 & 0 & 0\\
            0 & -\alpha & 0 & 0\\
            0 & 0 & \alpha & 0\\
            0 & 0 & 0 & -\alpha
        \end{pmatrix}\quad {\rm{or}}\quad
        A_1=
        \begin{pmatrix}
            \alpha & 0 & 0 & 0\\
            0 & \alpha & 0 & 0\\
            0 & 0 & \alpha & 0\\
            0 & 0 & 0 & -\alpha
        \end{pmatrix},
    \end{equation}
    where $\alpha\neq0$ and cylindrical in $\xi_2$;
     \begin{equation}\label{eq:shape4}
        A_1 =
        \begin{pmatrix}
            \alpha & 0 & 0 & 0\\
            0 & \frac{\beta}{\alpha} & 0 & 0\\
            0 & 0 & \pm\frac{\beta}{\alpha} & 0\\
            0 & 0 & 0 & \pm\alpha
        \end{pmatrix}
        ,\quad
        A_2 =
        \begin{pmatrix}
            0 & 0 & 0 & 0\\
            0 & p  & 0 & 0\\
            0 & 0 & q & 0\\
            0 & 0 & 0 & 0
        \end{pmatrix},
    \end{equation}
    where $pq=\pm\left(\alpha^2  - \frac{\beta^2}{\alpha^2}\right)$ and $\alpha\beta\neq0$;
      \begin{equation}\label{eq:shape5}
        A_1 =
        \begin{pmatrix}
            \alpha & 0 & 0 & 0\\
            0 & \frac{\beta}{\alpha} & 0 & 0\\
            0 & 0 & \frac{\gamma}{\alpha} & 0\\
            0 & 0 & 0 & \frac{\delta}{\alpha}
        \end{pmatrix}
        ,\quad
        A_2 =
        \begin{pmatrix}
            0 & 0 & 0 & 0\\
            0 & p  & 0 & 0\\
            0 & 0 & q & 0\\
            0 & 0 & 0 & r
        \end{pmatrix},
    \end{equation}
    where $\alpha\neq 0$ and
    \begin{equation*}
        pq = \pm\delta - \frac{\beta\gamma}{\alpha^2},\quad qr = \pm\beta - \frac{\gamma\delta}{\alpha^2},\quad pr = \pm\gamma - \frac{\beta\delta}{\alpha^2}.
    \end{equation*}
\end{enumerate}
\begin{proof}
    Since the normal connection is flat, at each point $x\in M$, there exists some orthonormal basis $\{e_1,e_2,e_3,e_4\}$ such that the shape operators $A_1$ and $A_2$ are simultaneously diagonalizable. We fix the point $x\in M$. If $M$ is totally geodesic, then it is cylindrical. If $M$ is not totally geodesic, we can assume that $h(e_1,e_1)\neq0$. Then we put $\xi_1 = \frac{h(e_1,e_1)}{||h(e_1,e_1)||}$ and $\xi_2$ is the unit normal vector perpendicular to $\xi_1$. Then the shape operators $A_1$ and $A_2$ are respectively given by
    \begin{equation*}
        A_1 = \begin{pmatrix}
            a_1 & 0 & 0 & 0\\
            0 & a_2 & 0 & 0\\
            0 & 0 & a_3 & 0\\
            0 & 0 & 0 & a_4
        \end{pmatrix} (a_1\neq 0)
        ,\quad
        A_2 = \begin{pmatrix}
            0 & 0 & 0 & 0\\
            0 & b_2 & 0 & 0\\
            0 & 0 & b_3 & 0\\
            0 & 0 & 0 & b_4
        \end{pmatrix}
    \end{equation*}
    Using Lemma~\ref{lem:gST} and Gauss equation, we have that $M$ is weakly Einstein if and only if
    \begin{align}
        &a_1 a_2 = a,\label{eq:gst1}\\
        &a_1 a_3 = b,\label{eq:gst2}\\
        &a_1 a_4 = c,\label{eq:gst3}\\
        &a_3 a_4 + b_3 b_4 = a',\label{eq:gst4}\\
        &a_2 a_4 + b_2 b_4 = b',\label{eq:gst5}\\
        &a_2 a_3 + b_2 b_3 = c',\label{eq:gst6}
    \end{align}
    where $a=R_{1221}$, $a'=R_{3443}$, $b=R_{1331}$, $b'=R_{2442}$, $c=R_{1441}$, $c'=R_{2332}$, and $a^2 = a'^2$, $b^2 = b'^2$, $c^2 = c'^2$.

    For resolving this system, we first compute $a_1$: by using the equations~\eqref{eq:gst1}-\eqref{eq:gst3} and the equality $||R||^2 = 8(a^2 + b^2 + c^2)=8a_1^2 (a_2^2 +a_3^2+a_4^2)$, we can find that $a_1$ is a solution of the following equation:
    \begin{equation*}
        x^4 - \Tr\big((A_1)^2\big) x^2 + \frac{||R||^2}{8}=0.
    \end{equation*}
    Indeed, we observe that the discriminant $D$ of this equation (set $t:=x^2$) is
    \begin{equation*}
        D = \big\{\Tr\big((A_1)^2\big)\big\}^2 - \frac{||R||^2}{2}=(a_1^2 - a_2^2 - a_3^2-a_4^2)^2 \ge 0.
    \end{equation*}

    To determine other unknowns, we consider the rank of the Riemannian curvature operator at $x$, denoted by $\rank\mathcal{R}(x)$. Here, the Riemannian curvature operator at $x$ is a symmetric map $\mathcal{R}:\Lambda^2 T_x M \to \Lambda^2 T_x M$ defined by
    $\<\mathcal{R}(X\wedge Y),Z\wedge W\>=-R(X,Y,Z,W)$.

    \smallskip
    \noindent\textsc{Case 1}: $\rank\mathcal{R}(x)=0$. Then $M$ is flat at $x$. By \eqref{eq:gst1}-\eqref{eq:gst6}, $M$ is cylindrical.

    \smallskip
    \noindent\textsc{Case 2}: $\rank\mathcal{R}(x)=2$. Let us consider the case of $a$ and $a'$ are not equal to zero. Since $a_1\neq 0$, from \eqref{eq:gst2} and \eqref{eq:gst3} we obtain $a_3=a_4=0$. Then Equations~\eqref{eq:gst5} and \eqref{eq:gst6} imply that $b_2b_3=b_2b_4=0$. If $b_2 \neq 0$, then $b_3=b_4=0$, which implies $a'=0$ by \eqref{eq:gst4}, a contradiction. Hence $b_2 =0$. Therefore, the shape operators are of the form
    \begin{equation*}
        A_1 =
        \begin{pmatrix}
            \alpha & 0 & 0 & 0\\
            0 & \beta & 0 & 0\\
            0 & 0 & 0 & 0\\
            0 & 0 & 0 & 0
        \end{pmatrix},\quad
        A_2 =
        \begin{pmatrix}
            0 & 0 & 0 & 0\\
            0 & 0  & 0 & 0\\
            0 & 0 & p & 0\\
            0 & 0 & 0 & q
        \end{pmatrix},
    \end{equation*}
    where $(\alpha\beta)^2 = (p q)^2\neq 0$. For the remaining cases, we can apply similar arguments and obtain the same result by relabelling.

    \smallskip
    \noindent\textsc{Case 3}: $\rank\mathcal{R}(x)=4$. Consider the case $a=a'=0$. Then, from \eqref{eq:gst1} and \eqref{eq:gst4}, we get $a_2=0$ and $a_3a_4+b_3b_4=0$. Using \eqref{eq:gst2}, \eqref{eq:gst3}, \eqref{eq:gst5}, and \eqref{eq:gst6}, we obtain
    $(a_1a_3)^2 = (b_2b_4)^2$ and $(a_1a_4)^2=(b_2b_3)^2$. It follows that $a_1^4 a_3^2 a_4^2 = b_2^4 b_3^2 b_4^2 = b_2^4 a_3^2 a_4^2$, which implies that $a_1^2 = b_2^2$ since $a_3,a_4 \neq 0$. Hence, we have the form of the shape operators given by
    \begin{equation*}
        A_1 =
        \begin{pmatrix}
            \alpha & 0 & 0 & 0\\
            0 & 0 & 0 & 0\\
            0 & 0 & \frac{\beta}{\alpha} & 0\\
            0 & 0 & 0 & \frac{\gamma}{\alpha}
        \end{pmatrix},\quad
        A_2 =
        \begin{pmatrix}
            0 & 0 & 0 & 0\\
            0 & \pm\alpha  & 0 & 0\\
            0 & 0 & p & 0\\
            0 & 0 & 0 & q
        \end{pmatrix},
    \end{equation*}
    where $pq = -\frac{\beta\gamma}{\alpha^2}\neq 0$. For the remaining cases, we can apply similar arguments and obtain the same result by relabelling.

    \smallskip
    \noindent\textsc{Case 4}: $\rank\mathcal{R}(x)=6$. In this situation, we divide the cases for the rank of $A_2$ (write $\rank A_2$).

    First, we assume that $\rank A_2$ is $0$ or $1$. Then we get
    \begin{align}
        (a_1a_2)^2 &= (a_3a_4)^2,\label{eq:squaregst1}\\
        (a_1a_3)^2 &= (a_2a_4)^2,\label{eq:squaregst2}\\
        (a_1a_4)^2 &= (a_2a_3)^2.\label{eq:squaregst3}
    \end{align}
    We note that $a_2,a_3,a_4\neq0$. Multiplying $a_3^2$ to \eqref{eq:squaregst1} and using \eqref{eq:squaregst3}, we obtain $a_1^4 a_4^2= a_3^4 a_4^2$, which leads to $a_1 = a_3$ or $a_1 = -a_3$. Consider $a_1=a_3$. Using again \eqref{eq:squaregst3} we get $a_2=a_4$ or $a_2 =-a_4$. Then \eqref{eq:squaregst2} implies that $a_1=a_2$ or $a_1=-a_2$. For the case of $a_1=-a_3$, the same result is obtained. Hence, either $M$ is umbilical with respect to $\xi_1$ and cylindrical in $\xi_2$, or the possible forms of $A_1$ are
    \begin{equation*}
        \begin{pmatrix}
            \alpha & 0 & 0 & 0\\
            0 & -\alpha & 0 & 0\\
            0 & 0 & \alpha & 0\\
            0 & 0 & 0 & -\alpha
        \end{pmatrix}\; {\text{or}}\;
        \begin{pmatrix}
            \alpha & 0 & 0 & 0\\
            0 & \alpha & 0 & 0\\
            0 & 0 & \alpha & 0\\
            0 & 0 & 0 & -\alpha
        \end{pmatrix} (\alpha\neq 0).
    \end{equation*}

    Second, we assume that $\rank A_2 =2$. Without loss of generality, we let $b_4=0$. Then we get \eqref{eq:squaregst1}, \eqref{eq:squaregst2}, and
\begin{equation}\label{eq:squaregst4}
            (a_1a_4)^2 = (a_2a_3+b_2b_3)^2.
\end{equation}
Multiplying \eqref{eq:squaregst1} and \eqref{eq:squaregst2}, we obtain $a_1^4a_2^2a_3^2=a_4^4a_2^2a_3^2$, which implies that $a_1=a_4$ or $a_1=-a_4$. Then, from \eqref{eq:squaregst1}, we have $a_2=a_3$ or $a_2=-a_3$. By \eqref{eq:squaregst4}, we have $a_1^4 = (\pm a_2^2 + b_2b_3)^2$. Therefore, the shape operators are of the form
 \begin{equation*}
        A_1 =
        \begin{pmatrix}
            \alpha & 0 & 0 & 0\\
            0 & \frac{\beta}{\alpha} & 0 & 0\\
            0 & 0 & \pm\frac{\beta}{\alpha} & 0\\
            0 & 0 & 0 & \pm\alpha
        \end{pmatrix}
        ,\quad
        A_2 =
        \begin{pmatrix}
            0 & 0 & 0 & 0\\
            0 & p  & 0 & 0\\
            0 & 0 & q & 0\\
            0 & 0 & 0 & 0
        \end{pmatrix},
    \end{equation*}
    where $pq=\pm\left(\alpha^2  - \frac{\beta^2}{\alpha^2}\right)$ and $\alpha\beta\neq0$.

    Third, we assume that $\rank A_2 =3$. Then, the form of the shape operators is the following:
    \begin{equation*}
        A_1 =
        \begin{pmatrix}
            \alpha & 0 & 0 & 0\\
            0 & \frac{\beta}{\alpha} & 0 & 0\\
            0 & 0 & \frac{\gamma}{\alpha} & 0\\
            0 & 0 & 0 & \frac{\delta}{\alpha}
        \end{pmatrix}
        ,\quad
        A_2 =
        \begin{pmatrix}
            0 & 0 & 0 & 0\\
            0 & p  & 0 & 0\\
            0 & 0 & q & 0\\
            0 & 0 & 0 & r
        \end{pmatrix},
    \end{equation*}
    where $\alpha\neq 0$ and
\begin{equation*}
        pq = \pm\delta - \frac{\beta\gamma}{\alpha^2},\quad qr = \pm\beta - \frac{\gamma\delta}{\alpha^2},\quad pr = \pm\gamma - \frac{\beta\delta}{\alpha^2}.
    \end{equation*}
\end{proof}
\end{thm}

\begin{rem}
    Let $M=M_1(c)\times M_2(-c)$ be the product of two surfaces of constant Gaussian curvature $c$ and $-c$, where $c\neq 0$. It can be easily checked that $M$ is weakly Einstein and $M$ is isometrically immersed in $\mathbb{R}^6$ whose normal connection is flat. The shape operators correspond to the case~\eqref{eq:shape1}.
\end{rem}

\section*{Acknowledgements}
This work was supported by the National Research Foundation of Korea (NRF) grant funded by the Korea government (MSIT) (NRF-2019R1A2C1083957). The authors thank Prof. Yuri Nikolayevsky for several valuable comments and careful proofreading of this paper.



\end{document}